# Some improved Gaussian correlation inequalities for symmetrical n-rectangles extended to some multivariate gamma distributions and some further probability inequalities

Thomas Royen

TH Bingen, University of Applied Sciences,

Berlinstrasse 109, D55411 Bingen, Germany

e-mail: thomas.royen@t-online.de

**Abstract.** The Gaussian correlation inequality (GCI) for symmetrical n-rectangles is improved if the absolute components have a joint cumulative distribution function (cdf), which is MTP$_2$ (multivariate totally positive of order 2). Inequalities of the here given type hold at least for all MTP$_2$ - cdfs on $\mathbb{R}^n$ or $\mathbb{R}^n_+$ with everywhere positive smooth densities. In particular, at least some infinitely divisible multivariate chi-square distributions (gamma distributions in the sense of Krishnamoorthy and Parthasarathy) with any positive real "degree of freedom" are shown to be MTP$_2$. Moreover, further numerically calculable probability inequalities for a broad class of multivariate gamma distributions are derived and a different improvement for inequalities of the GCI-type - and of a similar type with three instead of two groups of components - with more special correlation structures. The main idea behind these inequalities is to find for a given correlation matrix with positive correlations a further correlation matrix with smaller correlations whose inverse is an M-matrix and where the corresponding multivariate gamma distribution function is numerically available.

## 1. Introduction

Let $(X_1,...,X_n)$ be a random vector on $\mathbb{R}^n$ or $\mathbb{R}^n_+ = (0,\infty)^n$ with a continuous everywhere positive probability density function (pdf) and a cumulative distribution function (cdf)

$$F(x_1,...,x_n) = P\left(\cap_{i=1}^n A_i\right),\ A_i = \{X_i \le x_i\}$$

which satisfies the inequality

$$P(\cap_{i=1}^n A_i) \ge P(\cap_{i\in I} A_i)P(\cap_{i\notin I} A_i),\ \varnothing \ne I \subset \{1,...,n\}, \tag{1.1}$$

for all values $x_i$. Such a random vector is called here "completely positively lower orthant-dependent" (CPLOD), whereas "PLOD" means only the special case of (1.1) with $I$ of size 1 (or n-1) , which entails only the inequality

$$P(\cap_{i=1}^n A_i) \ge \prod_{i=1}^n P(A_i).$$

Alternatively, the CPLOD-property can also be denoted by MPQD, i.e. "multivariate positively quadrant dependent".

The inequality (1.1) is in particular the Gaussian correlation inequality for symmetrical n-rectangles (GCI) if the $X_i$ are the absolute values of the components of a centered normal random vector. This inequality was – after some earlier essential steps – completely proved in [21]. This GCI implies the general form of the GCI for centrally symmetrical convex regions, but this paper deals only with inequalities for rectangular regions. For an MTP2 - density ("multivariate totally positive of order 2", for the definition see section 3) the inequality (1.1) is always satisfied since (1.1) holds for "positively associated" random variables and since MTP2 implies positive association. However, a little bit more can be said if the cdf of $(X_1,...,X_n)$ is MTP2 . In section 3 a very simple proof is given for the following inequality for right-truncated $X_i$:

$$\frac{P(\bigcap_{i=1}^n A_i)}{P(\bigcap_{i\in I} B_i)} \geq \frac{P(\bigcap_{i\in I} A_i)}{P(\bigcap_{i\in I} B_i)} \frac{P(\bigcap_{i\notin I} A_i)}{P(\bigcap_{i\notin I} B_i)} \Leftrightarrow \frac{P(\bigcap_{i=1}^n A_i)}{P(\bigcap_{i\in I} A_i)P(\bigcap_{i\notin I} A_i)} \geq \frac{P(\bigcap_{i=1}^n B_i)}{P(\bigcap_{i\in I} B_i)P(\bigcap_{i\notin I} B_i)} \tag{1.2}$$

with $A_i = \{X_i \leq x_i\}$, $B_i = \{X_i \leq b_i\}$, $x_i \leq b_i$, $\varnothing \neq I \subset \{1,...,n\}$. This is an improvement of the inequality (1.1) within the region $B_1 \times ... \times B_n$ if the last fraction in (1.2) is greater than 1. E.g. for the GCI with a non-singular irreducible covariance matrix $\Sigma = (\sigma_{ij})$ (i.e. rank$(\sigma_{ij})_{i\in I, j\notin I}$ is always positive) this holds for all positive numbers $b_i$.

The GCI can be equivalently considered as an inequality for a chi-square or gamma random vector $(X_1,...,X_n)$ with $X_i = \frac{1}{2}Z_i^2$, where $(Z_1,...,Z_n)$ has an $N_n(0,\Sigma)$-distribution. The multivariate gamma distribution in the sense of Krishnamoorthy and Parthasarathy [13] coincides, apart from a scale factor 2, with the distribution of the diagonal of a Wishart $W_n(\nu,\Sigma)$- matrix with the "degree of freedom" $\nu = 2\alpha \in \mathbb{N}$ or $\nu > n-1$. This $\Gamma_n(\alpha,\Sigma)$ - distribution has the Laplace transform (Lt)

$$\left|I_n + \Sigma T\right|^{-\alpha} \tag{1.3}$$

with the $(n \times n)$- identity matrix $I_n$, the "associated" covariance matrix $\Sigma$ (not to confuse with the covariance matrix of the components $X_i$) and $T = Diag(t_1,...,t_n)$, $t_i \geq 0$. In this paper $\Sigma$ is always supposed to be non-singular. Formula (1.3) provides too the Lt of a pdf ($\Gamma_n(\alpha,\Sigma)$ - density per definition) at least for all non-integer values $2\alpha > \left[(n-1)/2\right]$ (see [24]).

In [21] the inequality (1.1) for the absolute Gaussian random variables was extended to $\Gamma_n(\alpha,\Sigma)$ - distributions with $2\alpha \in \mathbb{N} \cup (n-2,\infty)$ and the same proof holds now for all non-integer values $2\alpha > \left[(n-1)/2\right]$ too.

For some special structures of $\Sigma$ the $\Gamma_n(\alpha,\Sigma)$ - distribution exists for smaller values of $\alpha$ (see (2.9) in section 2). All $\alpha > 0$ are admissible iff there exists a signature matrix $S = Diag(s_1,...,s_n)$, $s_i = \pm 1$, for which $S\Sigma^{-1}S$ is an M-matrix (see [1] and for an equivalent condition [7]), i.e. all the off-diagonal elements of $S\Sigma^{-1}S$ are non-positive and all elements of $S\Sigma S$ are non-negative. The Lt in (1.3) shows the invariance of the $\Gamma_n(\alpha,\Sigma)$ - distribution under the transformation $\Sigma \to S\Sigma S$.

If $(Z_1,...,Z_n)$ has an $\mathrm{N}_n(0,\Sigma)$-distribution, then $(|Z_1|,...,|Z_n|)$ is MTP2 iff there exists a signature matrix $S$, which generates an M-matrix $S\Sigma^{-1}S$, (see [11]), and MTP2 for $(|Z_1|,...,|Z_n|)$ is equivalent to MTP2 for the $\Gamma_n(\frac{1}{2},\Sigma)$-distribution. Therefore, we have the remarkable fact, that a $\Gamma_n(\frac{1}{2},\Sigma)$-distribution is infinitely divisible (inf. div.) iff it is MTP2. This might seduce to the (very optimistic) conjecture that all inf. div. $\Gamma_n(\alpha,\Sigma)$ - distributions are MTP2 too (or at least positively associated). The inf. div. gamma distribution with the density from (2.12) is shown to be MTP2 in section 4, where also two further examples with inf. div. MTP2 $\Gamma_n(\alpha,\Sigma)$ - distributions are given.. On the other hand, the inequality (1.1) for general inf. div. $\Gamma_n(\alpha,\Sigma)$ - distributions with non-integer values $2\alpha$ is proved until now for $2\alpha > [(n-1)/2]$, but for $2\alpha \in (0,\ [(n-1)/2)])$ only for an index set $I \subset \{1,...,n\}$ of size $1$ (or $(n-1)$), which entails only PLOD (see [22], [24]).

In the proof of (1.1) in [21] for a $\Gamma_n(\alpha,\Sigma)$ -distribution with $\Sigma=\begin{pmatrix}\Sigma_{11} & \Sigma_{12}\\ \Sigma_{21} & \Sigma_{22}\end{pmatrix}$ it is supposed that $\left|I_n+\Sigma_\vartheta T\right|^{-(\alpha+1)}$ is the Lt of a $\Gamma_n(\alpha+1,\Sigma_\vartheta)$ - pdf for all $\Sigma=\begin{pmatrix}\Sigma_{11} & \vartheta\Sigma_{12}\\ \vartheta\Sigma_{12} & \Sigma_{22}\end{pmatrix}$, $\vartheta\in(0,1)$. Furthermore, it was proved that the $\Gamma_n(\alpha,\Sigma_\vartheta)$ - cdf is strictly increasing in $\vartheta\in(0,1)$ for all positive fixed numbers $x_1,\ldots,x_n$ if $\operatorname{rank}(\Sigma_{12})>0$. But if $\Sigma^{-1}$ is an M-matrix, then frequently $\Sigma_\theta^{-1}$ is no M-matrix for some $\vartheta\in(0,1)$, and then, it is not sure if $\left|I_n+\Sigma_\vartheta T\right|^{-(\alpha+1)}$ is the Lt of a pdf with a non-integer value $2(\alpha+1)<\left[(n-1)/2\right]$, $n>6$. All these difficulties would disappear if $\left|I_n+\Sigma T\right|^{-\alpha}$ would be the Lt of a pdf for all $\alpha>1$ and all $\Sigma$.

In section 5, some different improvements of (1.1) are derived for some $\Gamma_n(\alpha,\Sigma)$ - distributions by means of inf. div. gamma distributions and the here cited theorem 4 from [22], which provides an essential generalization of a theorem of Bølviken and Joag-Dev in [3] for absolute normal Gaussian random vectors. In particular, some improvements for the corresponding GCI are obtained with $\alpha=\frac{1}{2}$. In some cases, the given positive lower bound for the "excess" $P(\cap_{i=1}^n A_i)-P(\cap_{i\in I}A_i)P(\cap_{i\notin I}A_i)$ is numerically available.

Furthermore, an approximation for the excess is recommended, which is in particular useful for small "p–values" $p=1-P\{\max X_i\le x\}$. This can be applied to many simultaneous or sequentially rejective statistical tests with multivariate chi-square test-statistics.

Formulas from the handbook of mathematical functions are cited by HMF and their numbers. $A>B$ for $(n\times n)$- matrices $(a_{ij})$, $(b_{ij})$ means here $a_{ij}\ge b_{ij}$ for all $i,j$ and $A\ne B$. $\Sigma_{(k)}$ means a summation over all decompositions $k=k_1+\ldots+k_p$ with $k_j\in\mathbb{N}_0$.

## 2. Some formulas for multivariate gamma distributions

Here some formulas are compiled, which are used in the subsequent sections. The most comprehensive collection for integral and series representations and approximations for $\Gamma_n(\alpha,R)$ - cdfs is possibly found in the appendix of [5], (see also [14],…,[20]). For non-central multivariate gamma distributions see [23]. Since the scaling is irrelevant for the aims of this paper, all the following formulas for $\Gamma_n(\alpha,\Sigma)$ - distributions are given for the standardized case, i.e. $\Sigma$ is now a non-singular correlation matrix $R=(r_{ij})$. The required scale factors can be easily inserted in applications with corresponding $\chi_n^2(2\alpha,\Sigma)$ - distributions.

The probability density $x^{\alpha-1}\mathrm{e}^{-x}/\Gamma(\alpha)$, $\alpha>0$, $x>0$, of a univariate standard gamma distribution is denoted by $g_\alpha(x)$ and the corresponding cdf by $G_\alpha(x)$. The non-central gamma-cdf with non-centrality parameter $y$ is given by

$$G_\alpha(x;y)=e^{-y}\sum\nolimits_{n=0}^{\infty}G_{\alpha+n}(x)\frac{y^n}{n!}\quad\text{with the density } g_\alpha(x;y). \tag{2.1}$$

For integer values $2\alpha$ we have

$$G_{1/2+n}(x;y)-\tfrac{1}{2}\left(erf(\sqrt{x}+\sqrt{y})+erf(\sqrt{x}-\sqrt{y})\right)=-\mathrm{e}^{-(x+y)}\sum\nolimits_{k=1}^{n}(xy^{-1})^{(k-1/2)/2}I_{k-1/2}(2\sqrt{xy}) \tag{2.2}$$
$$=-e^{-y}\sum\nolimits_{k=1}^{n}g_{1/2+k}(x)\,{}_0F_1(\tfrac{1}{2}+k;xy)$$

with the modified Bessel functions $I_{k-1/2}$, which are elementary functions (HMF 10.49.(ii)), and

$$G_{1+n}(x;y)-G_1(x;y)=-\mathrm{e}^{-y}\sum\nolimits_{k=1}^{n}g_{1+k}(x)\,{}_0F_1(1+k;xy) \tag{2.3}$$

or

$$G_n(x;y) = (\tfrac{x}{y})^{n/2}\,\frac{1}{\pi}\int_0^{\pi}\frac{y\cos(n\varphi)-\sqrt{xy}\cos((n-1)\varphi)}{x-2\sqrt{xy}\cos\varphi+y}\,e^{-(x-2\sqrt{xy}\cos\varphi+y)}d\varphi + G_0(x-y),\ n\in\mathbb{N}, \tag{2.4}$$

$$G_0(z) = \begin{Bmatrix} 0,\ z<0 \\ \frac{1}{2},\ z=0 \\ 1,\ z>0 \end{Bmatrix},$$

(see section 2 in [16]).

The functions $G_\alpha(x;y)$ can be extended to holomorphic functions of $y\in\mathbb{C}$. According to [4], the cdf $G_\alpha(x;y)$ is a strict log-concave function of $x$ for all $\alpha>1,\ y>0$ (for $0<\alpha\le 1$ see theorem 3.4 in [6]), and $G_\alpha(x;y)$ is strict TP2 for all $\alpha>0,$ which is equivalent to

$$\frac{\partial^2}{\partial y\partial x}\log G_\alpha(x;y)>0. \tag{2.5}$$

The 2nd fact is an easy consequence of the first one, because

$$\frac{\partial}{\partial x}\left(\frac{g_{\alpha+1}(x;y)}{G_{\alpha+1}(x;y)}+1\right)<0 \text{ implies}$$

$$0<\frac{\partial}{\partial x}\frac{G_{\alpha+1}(x;y)}{G_{\alpha+1}(x;y)+g_{\alpha+1}(x;y)}=\frac{\partial}{\partial x}\frac{G_{\alpha+1}(x;y)}{G_\alpha(x;y)}=\frac{\partial}{\partial x}\frac{G_{\alpha+1}(x;y)-G_\alpha(x;y)}{G_\alpha(x;y)}=\frac{\partial^2}{\partial x\partial y}\log G_\alpha(x;y)=$$

$$\frac{\partial}{\partial y}\frac{g_\alpha(x;y)}{G_\alpha(x;y)}. \tag{2.6}$$

The char. function of $g_\alpha(x;y)$ is $\hat{g}_\alpha(t;y)=(1-it)^{-\alpha}\exp(ity/(1-it)).$

From $\lim_{\varepsilon\to 0}\hat{g}_\alpha(\varepsilon t;\varepsilon^{-1}y)=\exp(ity)$ it follows

$$\lim_{\varepsilon\to 0}G_\alpha(\varepsilon^{-1}x;\varepsilon^{-1}y)=\begin{Bmatrix} 1,\ x>y \\ 0,\ x<y\end{Bmatrix}. \tag{2.7}$$

A correlation matrix $R_{n\times n}$ is called (positively) “m-factorial”, if m is the lowest integer allowing a representation

$$R = D + AA' \tag{2.8}$$

with a pos. def. $D=Diag(d_1,...,d_n)$ and an $(n\times m)$ - matrix $A$ of rank $m<n.$

With $D^{-1/2}RD^{-1/2}=I_n+BB'$ and the columns $\vec{b}_j$ of $B'$ we have for the $\Gamma_n(\alpha,R)$ - cdf the integral representation

$$F(x_1,...,x_n;\alpha,R)=E\left(\prod\nolimits_{j=1}^{n}G_\alpha(d_j^{-1}x_j;\tfrac{1}{2}\vec{b}_j'S\vec{b}_j)\right),\ 2\alpha\in\mathbb{N}\cup(m-1,\infty), \tag{2.9}$$

if $B$ is real, where the expectation refers to the random $W_m(2\alpha,I_m)$ - Wishart (or pseudo-Wishart) matrix $S,$ (see [16] or [18]). If there are one or more pure imaginary columns in $B,$ caused by some negative eigenvalues of $D^{-1/2}RD^{-1/2}-I_n,$ the $\Gamma_n(\alpha,R)$ - cdf is obtained at least for $2\alpha\in\mathbb{N}$ or $2\alpha>\max\left(\left[\frac{n-1}{2}\right],m-1\right).$ The verification of (2.9) uses the Lt of non-central univariate gamma densities, followed by integration over $S.$

The most simple case arises for a “one-factorial” correlation matrix

$$R=Diag(...,1-a_j^2,...)+\vec{a}\vec{a}' \tag{2.10}$$

with a real or imaginary column $\vec{a} = (a_1,...,a_n)'$ and $\max a_j^2 < 1$. Then

$$F(x_1,...,x_n;\alpha,R) = \int_0^\infty \Big(\prod\nolimits_{j=1}^n G_\alpha((1-a_j^2)^{-1}x_j;(1-a_j^2)^{-1}a_j^2 y)\Big) g_\alpha(y)dy,\ \ \alpha > 0, \tag{2.11}$$

if $\vec{a}$ is real. This cdf satisfies the inequality (1.1) and the corresponding inequality where all the $A_i$ are replaced by their complements $\overline{A}_i$ due to Kimball's inequality (see e.g. [25]). Because of (2.7) the limit case with $a_k^2 = \max a_k^2 = 1 > \max_{j \neq k} a_j^2$ is admissible. With an imaginary column $\vec{a}$ this distribution is not infinitely divisible, i.e. not all $\alpha > 0$ are admissible. For the extension to $\max a_j^2 > 1$ see [5] or [18].

Besides, there is a simple formula for the $\Gamma_n(\alpha,R)$ - pdf if $R^{-1} = (r^{ij})$ is of a "tree type". This means that the graph with the vertices $\{1,...,n\}$ and $n-1$ edges, exactly corresponding to the $n-1$ off-diagonal elements $r^{ij} \neq 0$ in $R^{-1}$, is a "spanning tree", i.e. connected without any cycles. The special case with a tridiagonal $R^{-1}$ and $2\alpha \in \mathbb{N}$ goes back to [2] and was generalized later in [15]. Such a $\Gamma_n(\alpha,R)$ - pdf is given by

$$f(x_1,...,x_n;\alpha,R) = \Big(|R| \prod\nolimits_{i=1}^n r^{ii}\Big)^{-\alpha} \Big(\prod\nolimits_{i=1}^n g_\alpha(r^{ii}x_i)\Big) \prod\nolimits_{i<j} {}_0F_1(\alpha; r^{ij2}x_i x_j). \tag{2.12}$$

## 3. Some properties of MTP2 – distributions

In this section the inequality (1.2) is shown for MTP2 - cdfs with a smooth density $f > 0$ on $\mathbb{R}^n$ or on $\mathbb{R}_+^n$. To make this paper more self-contained, some here used well known facts on MTP2 - functions are compiled with simple proofs, which are possibly of some interest for the non-specialized reader. For more theorems on MTP2 see e.g. [10] and [11].

Let $F(\vec{x}) = F(x_1,...,x_n)$ be an everywhere positive function on $D = D_1 \times ... \times D_n \subseteq \mathbb{R}^n$ with any intervals $D_i$. For any pair $(\vec{x},\vec{y}) \in D \times D$ the notation $\vec{x} \vee \vec{y}$ stands for $(\max(x_1,y_1),...,\max(x_n,y_n))$ and $\vec{x} \wedge \vec{y}$ for $(\min(x_1,y_1),...,\min(x_n,y_n))$. $F$ is per definition MTP2 on $D$ if

$F(\vec{x} \vee \vec{y})F(\vec{x} \wedge \vec{y}) \geq F(\vec{x})F(\vec{y})$ for all $(\vec{x},\vec{y}) \in D \times D.$

$F$ is MTP2 on $D$ iff $F$ is "TP2 in pairs", i.e. $F$ as a function of $(x_i,x_j) \in D_i \times D_j$, $i \neq j$, is always TP2 on $D_i \times D_j$ with any fixed values $x_k \in D_k$, $k \neq i,j$.

Now let $D$ be an open (bounded or unbounded) n-rectangle and $F \in C^2(D)$ an everywhere positive function. Then $F$ is MTP2 on $D$ iff

$$\frac{\partial^2}{\partial x_i \partial x_j} \log F(x_1,...,x_n) \geq 0 \text{ for all } i \neq j \text{ and all } \vec{x} \in D, \tag{3.1}$$

which is verified by integration over $x_i$ and $x_j$. A simple consequence is

**Theorem 1.** If the cdf $F \in C^2(\mathbb{R}^n)$ (or $F \in C^2(\mathbb{R}_+^n)$) is MTP2 with an everywhere positive density $f$ on $\mathbb{R}^n$ (or $f$ on $\mathbb{R}_+^n$), then $F$ satisfies the inequality (1.2).

**Proof.** The inequality (1.2) holds if $P(\cap_{k=1}^n A_k)(P(\cap_{k \in I} A_k)P(\cap_{k \notin I} A_k))^{-1}$ is non-increasing on $x_1 \leq b_1$,..., $x_n \leq b_n$ for every $I \subset \{1,...,n\}$. With the above assumptions this is equivalent to

$$\frac{\partial}{\partial x_i} \log P(\cap_{k=1}^n A_k) \leq \frac{\partial}{\partial x_i} \log P(\cap_{k \in I} A_k) \text{ for every } i \in I. \tag{3.2}$$

It is

$$\frac{\partial^2}{\partial x_j \partial x_i} \log F(x_1,...,x_n) \geq 0 \text{ for all pairs } (x_i, x_j),\ i \neq j,$$

and therefore

$$\frac{\partial}{\partial x_j}\left(\frac{\partial}{\partial x_i} \log P(\cap_{k=1}^n A_k)\right) \geq 0,\ i \in I, j \notin I.$$

Since $\lim_{\min x_j \to \infty, j \notin I} P(\cap_{k=1}^n A_k) = P(\cap_{k \in I} A_k)$,

the inequality (3.2) is satisfied and therefore the inequality (1.2) too. □

For three events $A = \cap_{k=1}^n A_k, B = \cap_{k=n_1+1}^{n_1+n_2} A_k, C = \cap_{k>n_1+n_2} A_k,\ n_1 + n_2 + n_3 = n,$ it follows e.g. with $\min x_k \to \infty,\ k > n_1 + n_2,$ the "cancelling inequality"

$$\frac{P(ABC)}{P(AC)P(B)} \geq \frac{P(AB)}{P(A)P(B)} \Leftrightarrow \frac{P(ABC)}{P(AB)P(C)} \geq \frac{P(AC)}{P(A)P(C)} \tag{3.3}$$

which is equivalent to the "conditional GCI"

$$P(BC \mid A) \geq P(B \mid A)P(C \mid A). \tag{3.4}$$

Similar equations arise by permutations of $A, B, C.$ Further inequalities for these three blocks of events are e.g.

$$P(ABC) + P(A)P(B)P(C) \geq P(AB)P(C) + P(AC)P(B), \tag{3.5}$$

since $P(AC)P(B)(PABC))^{-1} \leq P(A)P(B)(P(AB))^{-1}$ implies

$$P(ABC) - P(AB)P(C) - (P(AC)P(B) - P(A)P(B)P(C)) =$$

$$P(ABC) - P(AB)P(C) - \left(P(ABC)P(AC)P(B)(P(ABC))^{-1} - P(AB)P(C)P(A)P(B)(P(AB))^{-1}\right) \geq$$

$$\left(P(ABC) - P(AB)P(C)\right)\left(1 - P(AC)P(B)(P(ABC))^{-1}\right) \geq 0.$$

The $\text{MPT}_2$ – condition for the cdf is sufficient for the inequalities in (3.3), (3.4) and (3.5) but not necessary. On the other hand these inequalities do not hold for the example with $A = \{Z_1^2 \leq 5.2\}, B = \{Z_2^2 \leq 5.6], C = \{Z_3^2 \leq 5.6\}$ where $(Z_1, Z_2, Z_3)$ is an $N(0,R)$-vector with the correlations $r_{12} = \frac{3}{4}, r_{13} = \frac{1}{2}, r_{23} = 0.$

If $f(x_1,...,x_n)$ is an $\text{MTP}_2$ - pdf on $D$ then all m-dimensional marginal densities, $2 \leq m \leq n-1,$ are $\text{MTP}_2$ too (see e.g. (1.16) in [10]). A short proof is also found in chapter 2 in [9]. (Obviously, $f$ can be replaced by $cf$ with any constant $c > 0.$) If $f$ is a continuous $\text{MTP}_2$ - pdf on $D$ then the corresponding cdf $F$ is $\text{MTP}_2$ on $D$ too. This is a direct consequence of the following known fact, which is established here formally as a theorem. It is proved here for simplicity only for the important cases $D = \mathbb{R}^n$ or $D = \mathbb{R}_+^n,$ but the proof can be extended to different n-rectangles $D$ too. For more general versions see also [12].

Let $F \in \mathrm{C}^2(\mathbb{R}^n)$ or $F \in \mathrm{C}^2(\mathbb{R}_+^n)$ be a mixture-cdf of one of the following two types:

(a) $F(x_1,...,x_n) = \int_D \left(\prod_{k=1}^n G_k(x_k; y)\right) g(y)dy,$

where the $G_k$ are here any cdfs with everywhere positive pdfs $g_k$ on $\mathbb{R}$ for all $k$, or $g_k$ only positive on $\mathbb{R}_+$ for all $k$, additionally depending on a random parameter $Y$ with the density $g(y)$ on any interval $D \subset \mathbb{R}$.

(b) $F(x_1,...,x_n) = \int_D \left(\prod_{k=1}^n G_k(x_k; y_k)\right) g(y_1,...,y_n) dy_1...dy_n$, with any cdfs as before, but now depending on the values $y_k$ of a random parameter-vector $(Y_1,..,Y_n)$ with an MTP2 - pdf $g$ on any rectangular region $D = D_1 \times ... \times D_n$.

**Theorem 2.** The cdf $F$ from (a) is MTP2 on $\mathbb{R}^n$ (or $\mathbb{R}^n_+$) if all the derivatives $\frac{\partial}{\partial y} \frac{g_k(x;y)}{G_k(x;y)}$ are simultaneously not negative or simultaneously not positive for all $x$ and for all $y \in D$.

The cdf from (b) is MTP2 on $\mathbb{R}^n$ (or $\mathbb{R}^n_+$) if all the derivatives $\frac{\partial}{\partial y_k} \frac{g_k(x_k;y_k)}{G_k(x_k;y_k)}$ are simultaneously not negative or simultaneously not positive for all $x$ and for all $y_k \in D_k$.

**Proof.** (a) With a change of the notations $y, z$ in the following double integral we find

$$F\frac{\partial^2 F}{\partial x_j \partial x_i} - \frac{\partial F}{\partial x_i}\frac{\partial F}{\partial x_j} = \iint_{D\times D} \Big(\frac{g_j(x_j;y)}{G_j(x_j;y)} - \frac{g_j(x_j;z)}{G_j(x_j;z)}\Big) \frac{g_i(x_i;y)}{G_i(x_i;y)} \left(\prod_{k=1}^n G_k(x_k;y) G_k(x_k;z)\right) g(y)g(z)\, dydz =$$

$$\frac{1}{2}\iint_{D\times D} \Big(\frac{g_j(x_j;y)}{G_j(x_j;y)} - \frac{g_j(x_j;z)}{G_j(x_j;z)}\Big)\Big(\frac{g_i(x_i;y)}{G_i(x_i;y)} - \frac{g_i(x_i;z)}{G_i(x_i;z)}\Big) \left(\prod_{k=1}^n G_k(x_k;y) G_k(x_k;z)\right) g(y)g(z) dydz \geq 0,$$

and the inequality is strict if the integrand does not vanish identically. $F$ is MTP2 follows now from (3.1).

(b) Only for simplicity of notation the inequality $F\frac{\partial^2 F}{\partial x_j \partial x_i} - \frac{\partial F}{\partial x_i}\frac{\partial F}{\partial x_j} \geq 0$ is shown here only for $i = 1$ and $j = 2$.

For the left side we obtain

$$\iint_{D\times D} \Big(\frac{g_2(x_2;y_2)}{G_2(x_2;y_2)} - \frac{g_2(x_2;z_2)}{G_2(x_2;z_2)}\Big) \frac{g_1(x_1;y_1)}{G_1(x_1;y_1)} \prod_{k=1}^n G_k(x_k;y_k) G_k(x_k;z_k) \times$$

$$g(y_1,...,y_n) g(z_1,...,z_n) dy_1...,dy_n dz_1...,dz_n =$$

$$\frac{1}{2}\int_{D_1\times D_2}\int_{D_1\times D_2} \Big(\frac{g_1(x_1;y_1)}{G_1(x_1;y_1)} - \frac{g_1(x_1;z_1)}{G_1(x_1;z_1)}\Big)\Big(\frac{g_2(x_2;y_2)}{G_2(x_2;y_2)} - \frac{g_2(x_2;z_2)}{G_2(x_2;z_2)}\Big) h(y_1,y_2) h(z_1,z_2)\, dy_1 dy_2 dz_1 dz_2, \quad (3.6)$$

$$h(y_1,y_2) = \int_{\times_{k=3}^n D_k} \prod_{k=1}^n G_k(x_k;y_k)\, g(y_1,...,y_n)\, dy_3...dy_n.$$

Since the following step holds for all fixed values $x_k$, the dependence of $h$ on $x_1,...,x_n$ is not explicitely noted. The function $h(y_1,y_2)$ is TP2 on $D_1 \times D_2$. Therefore, we can write the integral in (3.6) as an integral over the region $D^+_{1,2}$, defined by $(y_1 - z_1)(y_2 - z_2) > 0$. This leads, after changing the notation $(y_2, z_2)$ to $(z_2, y_2)$ to

$$\int_{D^+_{1,2}} \Big(\frac{g_1(x_1;y_1)}{G_1(x_1;y_1)} - \frac{g_1(x_1;z_1)}{G_1(x_1;z_1)}\Big)\Big(\frac{g_2(x_2;y_2)}{G_2(x_2;y_2)} - \frac{g_2(x_2;z_2)}{G_2(x_2;z_2)}\Big) \times$$

$$\big(h(y_1,y_2)h(z_1,z_2) - h(y_1,z_2)h(z_1,y_2)\big) dy_1 dy_2 dz_1 dz_2 \geq 0. \ \square$$

As a by-product it follows from (2.7) with the non-central gamma cdf $G_k(x_k;y_k)=G_\alpha(\varepsilon^{-1}x_k;\varepsilon^{-1}y),\ \varepsilon\to 0,$ and any continuous MTP$_2$- density $g(y_1,...,y_n)$ on $\mathbb{R}^n_+$ that the cdf $G$ with the density $g$ is MTP$_2$ on $\mathbb{R}^n_+$ too.

## 4. Some examples for infinitely divisible MTP$_2$ gamma distributions

All the three examples of inf. div. $\Gamma_n(\alpha,R)$ - cdfs in this section can be shown to be MTP$_2$ using theorem 2 or the more general theorem 2.1 in [12] and the proof below in the 2nd example, where ${}_0F_1(\alpha;cxy),\ \alpha,c>0,$ is shown to be TP$_2$ on $(x,y)\in(0,\infty)^2$.

Because of (2.5) the $\Gamma_n(\alpha,R)$ - cdf $F$ from (2.11) with $a_i\in(0,1),\ i=1,...,n,$ satisfies the conditions (a) from theorem 2. Therefore, $F$ is MTP$_2$ and satisfies the inequality (1.2). This can also be shown without to use the MTP$_2$ - property by a direct verification of the inequality (3.2). If $Y_1,Y_2$ are positive random variables with everywhere positive pdfs $f_1$ and $f_2=f_1w$ on $\mathbb{R}_+$ and $\frac{\partial}{\partial y}w(y)>0$ on $\mathbb{R}_+,$ then it follows for the corresponding cdfs $1-F_2(y)>1-F_1(y)$ for all $y>0,$ and we find for every increasing function $\psi(y)$ on $\mathbb{R}_+$ with $E(\psi(Y_i))<\infty,\ i=1,2,$

$$\int_0^\infty \psi(y)(f_2(y)-f_1(y))dy>0.$$

With the functions

$$G_k(x_k;y)=G_\alpha((1-a_k^2)^{-1}x_k;(1-a_k^2)^{-1}a_k^2y),$$

$$f_1(y)=\left(\int_0^\infty\left(\prod\nolimits_{k=1}^n G_k(x_k;y)\right)g_\alpha(y)dy\right)^{-1}\left(\prod\nolimits_{k=1}^n G_k(x_k;y)\right)g_\alpha(y),$$

and $f_2(y)$ as $f_1(y),$ but only with the indices $k\in I,$ $w(y)$ proportional to the increasing factor $\left(\prod_{k\notin I}G_k(x_k;y)\right)^{-1}$ with the suitable norming factor and $\psi(y)=\frac{\partial}{\partial x_i}\log G(x_i;y),\ i\in I,$ the inequality (3.2) is satisfied and consequently (1.2) too.

Now the $\Gamma_n(\alpha,R)$ - pdf from (2.12) is shown to be MTP$_2$. The function

$$F(\alpha,z)={}_0F_1(\alpha;z)/\Gamma(\alpha)=\sum\nolimits_{k=0}^\infty\frac{z^k}{\Gamma(\alpha+k)k!}$$

has the derivative $F'(\alpha;z)=F(\alpha+1;z).$ To prove $\frac{\partial^2}{\partial x_i\partial x_j}\log f(x_1,...,x_n;\alpha,R)\geq 0,\ i\neq j,$ it is sufficient to show that

$$\left(\frac{\partial^2}{\partial x_j\partial x_i}F(\alpha;cx_ix_j)\right)F(\alpha;cx_ix_j)-\frac{\partial}{\partial x_i}F(\alpha;cx_ix_j)\frac{\partial}{\partial x_j}F(\alpha;cx_ix_j)\geq 0,\ \text{if } c>0.$$

This is equivalent to

$$F(\alpha;z)F(\alpha+1;z)+z\left(F(\alpha+2;z)F(\alpha;z)-(F(\alpha+1;z))^2\right)\geq 0,\ z\geq 0.$$

With HMF 10.31.3 we find with the modified Bessel functions $I_\mu(z)=(\frac{1}{2}z)^\mu F(\mu+1;\frac{1}{4}z^2)$ the multiplication formula

$$F(\mu;z)F(\nu;z)=\sum\nolimits_{k=0}^\infty\frac{(\mu+\nu-1+k)_k}{\Gamma(\mu+k)\Gamma(\nu+k)}\frac{z^k}{k!},$$

and therefore

$F(\alpha;z)F(\alpha+1;z)+z\big(F(\alpha+2;z)F(\alpha;z)-(F(\alpha+1;z))^2\big)$ =

$$\sum\nolimits_{k=0}^{\infty}\frac{(2\alpha+k)_k}{\Gamma(\alpha+k)\Gamma(\alpha+1+k)}\,\frac{z^k}{k!}-\sum\nolimits_{k=0}^{\infty}\frac{(2\alpha+1+k)_k}{\Gamma(\alpha+k)\Gamma(\alpha+1+k)}\Big(\frac{1}{\alpha+k}-\frac{1}{\alpha+k+1}\Big)\,\frac{z^{k+1}}{k!}=$$

$$\frac{1}{\Gamma(\alpha)\Gamma(\alpha+1)}+\sum\nolimits_{k=1}^{\infty}\frac{(2\alpha+k)_{k-1}(2\alpha+k-1)}{\Gamma(\alpha+k)\Gamma(\alpha+1+k)}\,\frac{z^k}{k!}>0.$$

Thus, the corresponding $\Gamma_n(\alpha,R)$- cdf is MTP2 too and satisfies the inequality (1.2). More on total positivity of hypergeometric functions of $x_1x_2$ is found in [8].

The third example is a generalization of the first one. We consider a $\Gamma(\alpha,R)$- distribution with the correlation matrix

$$R=R_{\vartheta}=\begin{pmatrix} R_{11} & \vartheta R_{12}\\ \vartheta R_{21} & R_{22}\end{pmatrix}=D+\begin{pmatrix}\vec{a}_1\vec{a}_1{}' & \vartheta\vec{a}_1\vec{a}_2{}'\\ \vartheta\vec{a}_2\vec{a}_1{}' & \vec{a}_2\vec{a}_2{}'\end{pmatrix},\ \vec{a}_1'=(a_1,\ldots,a_{n_1}),\ \vec{a}_2'=(a_{n_1+1},\ldots,a_n), \tag{4.1}$$

$1<\mathrm{n}_1<n-1,\ a_i\in(0,1),$ diagonal $D=D_1\oplus D_2=Diag(\ldots,1-a_i^2,\ldots),\ 0<\vartheta<((1+q_1^{-1})(1+q_2^{-1}))^{1/2},$

$q_i=\vec{a}_i'D_i^{-1}\vec{a}_i,\ i=1,2.$

Then

$$R_{\vartheta}^{-1}=D^{-1}-\frac{1}{(1+q_1)(1+q_2)-\vartheta^2q_1q_2}\begin{pmatrix}(1+(1-\vartheta^2)q_2)\vec{b}_1\vec{b}_1' & \vartheta\vec{b}_1\vec{b}_2'\\ \vartheta\vec{b}_2\vec{b}_1' & (1+(1-\vartheta^2)q_1)\vec{b}_2\vec{b}_2{}'\end{pmatrix},\ \vec{b}_i=D_i^{-1}\vec{a}_i, \tag{4.2}$$

is an M-matrix if $\vartheta\in[0,(\min(1+q_1^{-1},1+q_2^{-1}))^{1/2})$, (see also sec. 4 in [11]). Therefore, with these values $\vartheta$ the $\Gamma_n(\alpha,R_{\vartheta})$- cdf exists for all $\alpha>0$ and the same is true for all the corresponding $\Gamma_n(\alpha,SR_{\vartheta}S)$- cdfs with any signature matrix $S$ or, with other words, for all the corresponding $\Gamma_n(\alpha,R_{\vartheta})$- cdfs with $|\vartheta|\le\min\big((1+q_1^{-1})^{1/2},\ (1+q_2^{-1})^{1/2}\big)$ and $a_i\in(-1,1)$. But by the following theorem these $\Gamma_n(\alpha,R_{\vartheta})$- cdfs are shown to be MTP2 if $|\vartheta|\le 1,$ and therefore they satisfy the inequality (1.2) too. It remains here an open question if this holds too for the admissible values $|\vartheta|>1.$ Furthermore, by means of these distributions and theorem 4 in [22], we shall obtain a different improvement of the inequality (1.1) for some correlation structures in section 5.

**Theorem 3**. The $\Gamma_n(\alpha,R_{\vartheta})$- cdf with the correlation matrix $R_{\vartheta}$ from (4.1) and $\alpha>0$ is given by

$$F(x_1,\ldots,x_n;\alpha,R_{\vartheta})=$$

$$\sum\nolimits_{k=0}^{\infty}\binom{\alpha+k-1}{k}^{-1}\vartheta^{2k}\int_0^{\infty}G_{(1)}(\vec{x}_1;y)L_k^{(\alpha-1)}(y)g_{\alpha}(y)dy\int_0^{\infty}G_{(2)}(\vec{x}_2;y)L_k^{(\alpha-1)}(y)g_{\alpha}(y)dy,\ |\vartheta|\le 1, \tag{4.3}$$

with $G_{(1)}(\vec{x}_1;y)=\prod\nolimits_{j\le n_1}G_{\alpha}((1-a_j^2)^{-1}x_j;(1-a_j^2)^{-1}a_j^2y),\ G_{(2)}(\vec{x}_2;y)$ likewise with $j>n_1,$ and the generalized Laguerre polynomials $L_k^{(\alpha-1)}$, or alternatively by

$$(1-\vartheta^2)^{\alpha}\int_{\mathbb{R}_+^2}G_{(1)}(\vec{x}_1;(1-\vartheta^2)y_1)G_{(2)}(\vec{x}_2;(1-\vartheta^2)y_2)\,{}_0F_1(\alpha;\vartheta^2y_1y_2)g_{\alpha}(y_1)g_{\alpha}(y_2)dy_1dy_2,\ |\vartheta|<1. \tag{4.4}$$

**Remarks**. The special case of formula (4.3) with identical $a_j^2 = r_1,\ j \le n_1,\ a_j^2 = r_2\ ,\ j > n_1,$ was already derived in [19]. The abs. convergence of the series in (4.3) with $|\vartheta| = 1$ follows from the orthogonal series

$\sum_{k=0}^{\infty} c_{i,k} L_k^{(\alpha-1)}(y)$ for $G_{(i)}(\vec{x}_i; y)$ with $c_{i,k} = \binom{\alpha+k-1}{k}^{-1} \int_0^\infty G_{(i)}(\vec{x}_i; y) L_k^{(\alpha-1)}(y) g_\alpha(y)dy$ and $\|G_{(i)}\|^2 = \sum_{k=0}^{\infty} \binom{\alpha+k-1}{k} c_{i,k}^2 \le 1,$ but for $|\vartheta| = 1$ we have the simpler formula (2.11).

**Proof.** With $Z = (I_n + DT)^{-1} = Diag(z_1, ..., z_n),\ z_j = (1 + (1 - a_j^2)t_j)^{-1},\ Z = Z_1 \oplus Z_2,\ T = T_1 \oplus T_2,$ and $\vec{c}_i = (T_i Z_i)^{1/2} \vec{a}_i$ we obtain for the determinant in the corresponding Lt

$$|I_n + R_\vartheta T| = |Z|^{-1} \left| I_n + \begin{pmatrix} \vec{c}_1 \vec{c}_1{}' & \vartheta \vec{c}_1 \vec{c}_2{}' \\ \vartheta \vec{c}_2 \vec{c}_1{}' & \vec{c}_2 \vec{c}_2{}' \end{pmatrix} \right|. \tag{4.5}$$

The matrix behind $I_n$ has rank 2. Therefore, it is

$|I_n + R_\vartheta T| = |Z|^{-1} (1 + q_1 + q_2 + (1 - \theta^2) q_1 q_2)$ with $q_i = \vec{c}_i{}' \vec{c}_i$ and consequently with $\delta_i = 1 + q_i$ :

$$|I_n + R_\vartheta T|^{-\alpha} = (\prod_{i=1}^{n} z_j^\alpha)(\delta_1 \delta_2)^{-\alpha} (1 - \vartheta^2 (1 - \delta_1^{-1})(1 - \delta_2^{-1}))^{-\alpha} =$$

$$(\prod_{i=1}^{n} z_j^\alpha) \sum_{k=0}^{\infty} \binom{\alpha+k-1}{k} \vartheta^{2k} \left( \sum_{m=0}^{k} (-1)^m \binom{k}{m} \delta_1^{-(\alpha+m)} \right) \left( \sum_{m=0}^{k} (-1)^m \binom{k}{m} \delta_2^{-(\alpha+m)} \right). \tag{4.6}$$

For the inversion of this Lt we write

$$(\prod_{j \le n_1} z_j^\alpha)\, \delta_1^{-(\alpha+m)} = (\prod_{j \le n_1} z_j^\alpha) \int_0^\infty \exp(-y \sum_{j \le n_1} a_j^2 t_j z_j) g_{\alpha+m}(y) dy =$$

$$(\prod_{j \le n_1} z_j^\alpha) \int_0^\infty \exp\Big(-y \sum_{j \le n_1} \frac{a_j^2 (1 - z_j)}{1 - a_j^2}\Big) g_{\alpha+m}(y) dy =$$

$$\int_0^\infty \prod_{j \le n_1} \exp\Big(-\frac{a_j^2 y}{1 - a_j^2}\Big) \left( \sum_{k=0}^{\infty} \Big(\frac{a_j^2 y}{1 - a_j^2}\Big)^k \frac{z_j^{\alpha+k}}{k!} \right) g_{\alpha+m}(y) dy.$$

Inversion, followed by integration over $x_1, ..., x_{n_1}$ yields

$\int_0^\infty \big(\prod_{j \le n_1} G_\alpha((1 - a_j^2)^{-1} x_j; (1 - a_j^2)^{-1} a_j^2 y)\big) g_{\alpha+m}(y) dy,$ and it is

$$\sum_{m=0}^{k} (-1)^m \binom{k}{m} g_{\alpha+m}(y) = \sum_{m=0}^{k} (-1)^m \binom{k}{m} \frac{y^m}{(\alpha)_m} g_\alpha(y) = \binom{\alpha+k-1}{k}^{-1} L_k^{(\alpha-1)}(y) g_\alpha(y).$$

The inversion of $(\prod_{j > n_1} z_j^\alpha)\, \delta_2^{-(\alpha+m)}$ is accomplished in the same way, and this proves formula (4.3).

From the Poisson kernel (see HMF 18.18.27)

$$K(y_1, y_2; \alpha, \vartheta) = \sum_{k=0}^{\infty} \binom{\alpha+k-1}{k}^{-1} \vartheta^{2k} L_k^{(\alpha-1)}(y_1) L_k^{(\alpha-1)}(y_2) =$$

$$\frac{\Gamma(\alpha)(\vartheta(\sqrt{y_1 y_2}))^{-(\alpha-1)}}{1 - \vartheta^2} \exp\Big(-\frac{\vartheta^2 (y_1 + y_2)}{1 - \vartheta^2}\Big) I_{\alpha-1}\Big(\frac{2\vartheta \sqrt{y_1 y_2}}{1 - \vartheta^2}\Big) = \tag{4.7}$$

$$(1 - \vartheta^2)^{-\alpha} \exp\Big(-\frac{\vartheta^2 (y_1 + y_2)}{1 - \vartheta^2}\Big)\, {}_0F_1\Big(\alpha; \frac{\vartheta^2 y_1 y_2}{(1 - \vartheta^2)^2}\Big), \quad |\vartheta| < 1,$$

the formula (4.4) follows after the substitutions $y_i \to (1-\vartheta^2)y_i,\ i=1,2.$ $\square$

In the 2nd example of this section it was already verified that

$\frac{\partial^2}{\partial y_2 \partial y_1}\log({}_0F_1(\alpha;\vartheta^2 y_1 y_2)) > 0,$ i.e. ${}_0F_1(\alpha;\vartheta^2 y_1 y_2)$ is TP2 on $\mathbb{R}^2_+$. Then the $\Gamma_n(\alpha,R_\vartheta)$ - cdf is MTP2 if $|\vartheta|<1$. This follows directly from a limit case of theorem 2.1(b) in [12].

For $\vartheta^2 \in (1,(1+q_1^{-1})(1+q_2^{-1}))$ no formula similar to (4.3) or (4.4) seems to be known, but we can use formula (2.9) with $m=2,$ since $\operatorname{rank}(R_\vartheta - D)=2.$ $D^{-1/2}R_\vartheta D^{-1/2}-I_n$ has the reduced char. polynomial $\lambda^2-(q_1+q_2)\lambda+(1-\vartheta^2)q_1q_2$ and the eigenvalues

$$\lambda_\pm = \frac{q_1+q_2}{2} \pm \left( \left(\frac{q_1-q_2}{2}\right)^2 + q_1 q_2 \vartheta^2 \right)^{1/2},\ \ q_i = q_i = \vec{a}_i' D_i^{-1}\vec{a}_i. \tag{4.8}$$

With $b_j = a_j d_j^{-1/2} = a_j(1-a_j^2)^{-1/2}$ (not $b_j = a_j d_j^{-1}$ as in (4.2)) and

$$\mathrm{v}_1 = (b_1,\dots b_{n_1})',\ \mathrm{v}_{2\pm} = \frac{\vartheta q_1}{\lambda_\pm - q_2}(b_{n_1+1},\dots,b_n)' = \frac{\lambda_\pm - q_1}{\vartheta q_2}(b_{n_1+1},\dots,b_n)'$$

the eigenvectors are $\vec{e}_\pm = \|\mathrm{v}_1 \oplus \mathrm{v}_{2\pm}\|^{-1}(\mathrm{v}_1 \oplus \mathrm{v}_{2\pm})$. This entails $D^{-1/2}R_\vartheta D^{-1/2}-I_n = BB'$ with the $(n\times 2)$ - matrix $B=(\sqrt{\lambda_+}\vec{e}_+,\sqrt{\lambda_-}\vec{e}_-)$.

If $S$ is a $W_2(2\alpha,I_2)$ - Wishart matrix, then

$$\tfrac{1}{2}S = \begin{pmatrix} Y_1 & \sqrt{Y_1Y_2}\cos\Phi \\ \sqrt{Y_1Y_2}\cos\Phi & Y_2 \end{pmatrix}$$

has the density

$$g_\alpha(y_1)g_\alpha(y_2)f(\varphi),\ f(\varphi) = (\sin^2\varphi)^{(\alpha-1)}(B(\tfrac{1}{2},\alpha-\tfrac{1}{2}))^{-1},\ 0<\varphi<\pi,\ \alpha>\tfrac{1}{2},$$

and we obtain from ( 2.9) with the columns $\vec{b}_j = (b_{j1},b_{j2})'$ from $B'$ the $\Gamma_n(\alpha,R_\vartheta)$ - cdf

$$F(x_1,\dots,x_n;\alpha,R_\vartheta) =$$
$$\int_{\mathbb{R}^2_+}\left(\int_0^\pi\left(\prod_{j=1}^n G_\alpha(d_j^{-1}x_j;b_{j1}^2y_1+b_{j2}^2y_2+2b_{j1}b_{j2}\sqrt{y_1y_2}\cos\varphi)\right)f(\varphi)d\varphi\right)g_\alpha(y_1)g_\alpha(y_2)dy_1dy_2. \tag{4.9}$$

The 2nd column in $B$ is imaginary since $\lambda_-<0$ if $|\vartheta|>1$. Therefore, the non-centrality parameters $\tfrac{1}{2}\vec{b}_j'S\vec{b}_j$ are complex if $|\vartheta|>1,$ but the integral over $\varphi$ is real because of the symmetry of $f(\varphi)$. In the limit case $\alpha=\tfrac{1}{2}$, the probability measure of the angle $\Phi$ becomes symmetrically concentrated at the points $0$ and $\pi$. This leads to

$$F(x_1,\dots,x_n;\tfrac{1}{2},R_\vartheta) = \int_{\mathbb{R}^2_+} Re\left(\prod_{j=1}^n G_{1/2}(d_j^{-1}x_j;b_{j1}^2y_1+b_{j2}^2y_2+2b_{j1}b_{j2}\sqrt{y_1y_2})\right)g_{1/2}(y_1)g_{1/2}(y_2)dy_1dy_2. \tag{4.10}$$

The non-central gamma functions $G_{1/2}$ with a complex non-centrality parameter can be computed again by the error function as in (2.2).

**5. A further improvement of the inequality (1.1) for some multivariate gamma distributions and some related inequalities and approximations**

With the here cited Theorem 4 from [22] some positive lower bounds for the excess in the inequality (1.1) can be derived for some $\Gamma_n(\alpha,R)$ - distributions. After a permutation of the random components we can choose the index set $I=\{1,...n_1\}$, $1<n_1<n-1$, in (1.1) and consider only correlation matrices $\Sigma=R$ since the scaling is irrelevant.

**Theorem 4**. Let $R=(r_{ij})$ and $R_0=(r_{0,ij})$ be two different non-singular $(n\times n)$ - correlation matrices with $r_{ij}\geq r_{0,ij}>0$ for all $i\neq j$ and $r_0^{ij}\leq 0$ for all off-diagonal elements of $R_0^{-1}$. Then the $\Gamma_n(\alpha,R_0+\tau(R-R_0))$ - cdfs are strictly increasing in $\tau\in[0,1]$ for all positive numbers $x_1,...,x_n$ and $2\alpha\in\mathbb{N}\cup\left(\left[(n-1)/2\right],\infty\right)$.
This implies

$$F(x_1,...,x_n;\alpha,R)-F(x_1,...,x_n;\alpha,R_0)>0. \tag{5.1}$$

**Remarks.** This holds too if there exists any signature matrix $S$ for which $SR_0S$ and $SRS$ satisfy the above conditions. A similar inequality was proved in [3] for absolute Gaussian random vectors (here included with $\alpha=\frac{1}{2}$) under the additional restriction that all the matrices $(R_0+\tau(R-R_0))^{-1}$ are M-matrices too, whereas this condition is here only required for $R_0^{-1}$.

As an immediate consequence of theorem 4 and (4.2) we obtain with the correlation matrices $R_\vartheta$ from (4.1), $\vartheta\in[0,\min(\sqrt{1+(\vec{a}_1'D_1^{-1}\vec{a}_1)^{-1}},\sqrt{1+(\vec{a}_2'D_2^{-1}\vec{a}_2)^{-1}})]$ and $R=(r_{ij})=D+\begin{pmatrix}\vec{a}_1\vec{a}_1' & R_{12}\\ R_{21} & \vec{a}_2\vec{a}_2'\end{pmatrix}\neq R_\vartheta$ with $r_{ij}\geq\vartheta a_i a_j$, $i\leq n_1$, $j>n_1$, the inequality

$$P_R(\cap_{i=1}^n A_i)-P(\cap_{i\leq n_1}A_i)P(\cap_{i>n_1}A_i)=F(x_1,...,x_n;\alpha,R)-F(x_1,...,x_n;\alpha,R_0)> \tag{5.2}$$

$$F(x_1,...,x_n;\alpha,R_\vartheta)-F(x_1,...,x_n;\alpha,R_0) \tag{5.3}$$

for all positive numbers $x_1,...,x_n$ and $2\alpha\in\mathbb{N}\cup\left(\left[(n-1)/2\right],\infty\right)$. It is not required that $R^{-1}$ is an M-matrix too. The here given lower bound in (5.3) for the excess in (5.2) is numerically available by the formula (2.11) if $\vartheta=1$, by the formulas (4.3) or (4.4) if $\vartheta\in(0,1)$ and by (4.9) if $\vartheta>1$.

Now we consider a decomposition of the components $X_\mu$ into $p$ successive blocks with the sizes $n_i$, $n_1+...+n_p=n$, and $(n\times n)$ - correlation matrices

$$R_\Theta=D+(\vartheta_{ij}\vec{a}_i\vec{a}_j') \tag{5.4}$$

with $(n_i\times n_j)$ - blocks $\vec{a}_i\vec{a}_j'$, diagonal $D$, all numbers $a_\mu\in(0,1)$ and a non-singular correlation matrix $\Theta=(\vartheta_{ij})$. The following Lemma helps to compute the Lt $\left|I_n+R_\Theta T\right|^{-\alpha}$ of the $\Gamma_n(\alpha,R_\Theta)$ - distribution.

**Lemma 1.** Let $(\vartheta_{ij}\vec{b}_i\vec{b}_j')$ be an $(n\times n)$ - matrix with $(n_i\times n_j)$ - blocks, $\vec{b}_i\in\mathbb{R}^{n_i}$, $i=1,...,p$, and $Q=Diag(q_1,...,q_p)$ with $q_i=\vec{b}_i'\vec{b}_i\neq 0$. Then the eigenvalues $\lambda_1,...,\lambda_p$ of $\Theta Q$ coincide with the non-zero eigenvalues of $(\vartheta_{ij}\vec{b}_i\vec{b}_j')$, and it follows for the determinants

$$|\,I_n+(\vartheta_{ij}\vec{b}_i\vec{b}_j')T\,|\;=\;|\,I_p+Q^{1/2}\Theta Q^{1/2}\,|. \tag{5.5}$$

**Proof.** The matrices $(\vartheta_{ij}\vec{b}_i\vec{b}_j')$ and $\Theta Q$ have the same rank $p$.

$$(\vartheta_{ij}\vec{b}_i\vec{b}_j')\begin{pmatrix} c_1\vec{b}_1 \\ \vdots \\ c_p\vec{b}_p \end{pmatrix} = \begin{pmatrix} (q_1\vartheta_{11}c_1 + \ldots + q_p\vartheta_{1p}c_p)\vec{b}_1 \\ \vdots \\ (q_1\vartheta_{p1}c_1 + \ldots + q_p\vartheta_{pp}c_p)\vec{b}_p \end{pmatrix} = \lambda \begin{pmatrix} c_1\vec{b}_1 \\ \vdots \\ c_p\vec{b}_p \end{pmatrix} \Leftrightarrow \Theta Q\vec{c} = \lambda\vec{c}. \ \square$$

The computation of $|\, I_n + R_\Theta T \,|$ is now reduced to the computation of a determinant of a $(p \times p)$ - matrix as follows: With notations similar as in (4.5), but now with $p \geq 2$, we obtain

$$|\, I_n + R_\Theta T \,| \ = \ |\, Z \,|^{-1} |\, I_n + (\vartheta_{ij}\vec{b}_i\vec{b}_j') \,|, \ \vec{b}_i = (T_i Z_i)^{1/2}\vec{a}_i,$$

and consequently

$|\, I_n + R_\Theta T \,| \ = \ |\, Z \,|^{-1} |\, I_p + Q^{1/2}\Theta Q^{1/2} \,|$ with the elements $q_i = \vec{a}_i' T_i Z_i \vec{a}_i$ in $Q$. With $\delta_i = 1 + q_i$ this leads to the Lt

$$|\, I_n + R_\Theta T \,|^{-\alpha} = \left(\prod\nolimits_{\mu=1}^{n} z_\mu^{\alpha}\right)\left(\prod\nolimits_{i=1}^{p} \delta_i^{-\alpha}\right) |\, \tilde{\Theta} \,|^{-\alpha}, \tag{5.6}$$

where $\tilde{\Theta}$ is the correlation matrix with the correlations $\tilde{\vartheta}_k = \tilde{\vartheta}_{ij} = \vartheta_{ij}\sqrt{(1-\delta_i^{-1})(1-\delta_j^{-1})}.$

In theorem 5 below the $\Gamma_n(\alpha, R_\Theta)$ - cdf is given for an $R_\Theta$ as in (5.4) with three blocks. It would be easy to write down the corresponding formulas for more blocks, based on the binomial series for $|\, \tilde{\Theta} \,|^{-\alpha}$, but the computing effort for the resulting multivariate series would be very high.

The binomial series for $|\, \tilde{\Theta} \,|^{-\alpha} = (1-\tilde{\vartheta})^{-\alpha}$ is abs. convergent, but if the char. function (cf) is used instead of the Lt, an additional condition is required for more than two blocks. The binomial series is used now with complex numbers with absolute values $\rho$. A simple sufficient condition for all $\rho < 1$ is

$$\sum_{1 \leq j < k \leq 3} \tilde{\vartheta}_{jk}^4 \leq 1, \tag{5.7}$$

now with the correlations

$$\tilde{\vartheta}_{jk} = \vartheta_{jk}\left(\sum\nolimits_{\mu \in I_j} \frac{a_\mu^2}{1-a_\mu^2}\right)^{1/2}\left(1+\sum\nolimits_{\mu \in I_j} \frac{a_\mu^2}{1-a_\mu^2}\right)^{-1/2}\left(\sum\nolimits_{\mu \in I_k} \frac{a_\mu^2}{1-a_\mu^2}\right)^{1/2}\left(1+\sum\nolimits_{\mu \in I_k} \frac{a_\mu^2}{1-a_\mu^2}\right)^{-1/2}, \tag{5.8}$$

which is derived in the appendix.

The main elements in the following series representation of the $\Gamma_n(\alpha, R_\Theta)$ - cdf are the coefficients

$$c_{j,k}(\vec{x}_j;\alpha,\vec{a}_j) = \binom{\alpha+k-1}{k}^{-1}\int_0^\infty \left(\prod\nolimits_{\mu \in I_j} G_\alpha((1-a_\mu^2)^{-1}x_\mu;(1-a_\mu^2)^{-1}a_\mu^2 y)\right)L_k^{(\alpha-1)}(y)g_\alpha(y)dy \tag{5.9}$$

of the orthogonal series with generalized Laguerre-polynomials $L_k^{(\alpha-1)}$ for the products of non-central gamma cdfs belonging to block $I_j$ .

**Theorem 5**. Let $R_\Theta$ be a correlation matrix with $p=3$ blocks as described in (5.4). With the notation $\vartheta_j = \vartheta_{ik}$ for the correlations of $\Theta$ the $\Gamma_n(\alpha, R_\Theta)$ - cdf is given by

$$F(x_1,\ldots,x_n;\alpha,R_\Theta) = \sum\nolimits_{k=0}^{\infty}\sum\nolimits_{(k)}(\alpha)_k \frac{(-2)^{k_4}}{k_4!}\prod\nolimits_{j=1}^{3}\frac{\vartheta_j^{2k_j+k_4}}{k_j!}c_{j,k-k_j}(\vec{x}_j;\alpha,\vec{a}_j) \tag{5.10}$$

under the condition (5.7).

**Remarks.** With $\vartheta_j \in (0,1),\ j=1,2,3,\ R_\Theta^{-1}$ is an M-matrix iff

$$(\vartheta_1 + (\vartheta_1 - \vartheta_2\vartheta_3)p_1)(\vartheta_2 + (\vartheta_2 - \vartheta_1\vartheta_3)p_2)(\vartheta_3 + (\vartheta_3 - \vartheta_1\vartheta_2)p_3) > 0 \tag{5.11}$$

with $p_j = \sum_{\mu \in I_j} \frac{a_\mu^2}{1-a_\mu^2}$ (see example 4.3 in [11]). Then all $\alpha > 0$ are admissible. A 2nd formula

$$F(x_1,\dots,x_n;\alpha,R_\Theta) = \prod_{j=1}^{3} c_{j,0}(\vec{x}_j;\alpha,\vec{a}_j) +$$

$$\sum_{N=2}^{\infty}(-1)^N \sum_{m_1+m_2+m_3=n-\varepsilon_N} \Big(\sum_{m=0}^{\min m_j} \frac{2^{2m+\varepsilon_N}(\alpha)_{n-m}}{(2m+\varepsilon_N)!\prod_{j=1}^{3}(m_j-m)!}\Big)\prod_{j=1}^{3} \vartheta_j^{2m_j+\varepsilon_N} c_{j,n-m_j}(\vec{x}_j;\alpha,\vec{a}_j) \tag{5.12}$$

with $n=[\frac{1}{2}N]$ and $\varepsilon_N = N \bmod 2$ is obtained by rearranging the series in (5.10) into a series of homogeneous polynomials $P_N(\vartheta_1,\vartheta_2,\theta_3)$ of degree $N,$ but for the convergence of this series a spectral radius $\|\tilde{\Theta} - I_3\| < 1$ is supposed.

**Proof of theorem 5**. It follows from (5.6), but now for the corresponding cf with complex values $\zeta_j = 1-\delta_j^{-1}$ and the condition (5.7)

$$|I_n - iR_\Theta T|^{-\alpha} = \Big(\prod_{\mu=1}^{n} z_\mu^{\alpha}\Big)\Big(\prod_{i=1}^{3}\delta_i^{-\alpha}\Big)\big(1-(\vartheta_1^2\zeta_2\zeta_3 + \vartheta_2^2\zeta_1\zeta_3 + \vartheta_3^2\zeta_1\zeta_2 - 2\vartheta_1\vartheta_2\vartheta_3\zeta_1\zeta_2\zeta_3)\big)^{-\alpha} =$$

$$\Big(\prod_{\mu=1}^{n} z_\mu^{\alpha}\Big)\Big(\prod_{i=1}^{3}\delta_i^{-\alpha}\Big)\sum_{k=0}^{\infty}\sum_{(k)}(\alpha)_k \frac{(-2)^{k_4}}{k_4!}\prod_{j=1}^{3}\frac{\vartheta_j^{2k_j+k_4}}{k_j!}\zeta_j^{k-k_j}. \tag{5.13}$$

Inversion and integration as in the proof of theorem 3 provides formula (5.10). □

If $R_\Theta$ is a correlation matrix as in (5.4) with $p=3$ blocks, which satisfies the condition (5.11), then it follows from theorem 4 for all correlation matrices $R > R_\Theta$ the inequality

$$P_R\big(\bigcap_{i=1}^{n} A_i\big) > P_{R_\Theta}\big(\bigcap_{i=1}^{n} A_i\big),\ 2\alpha \in \mathbb{N} \cup \big(\big[(n-1)/2\big],\infty\big),\ x_1,\dots,x_n > 0.$$

Suitable numbers $a_\mu \in (0,1)$ are sometimes found in the following way: At first the sums $\sum_{\mu<\nu}(\ln(a_\mu a_\nu r_{\mu\nu}^{-1}))^2$ are minimized for each index block $I_j$ of size $n_j \geq 3$ separately by

$$a_\mu = \left(\frac{\prod_{\nu \in I_j} r_{\mu\nu}}{\Big(\prod_{\mu<\nu} r_{\mu\nu}\Big)^{1/(n_j-1)}}\right)^{1/(n_j-2)}$$

(which is simpler than minimizing $\sum(r_{\mu\nu} - a_\mu a_\nu)^2$ ). Then the $a_\mu$ are replaced by $\lambda_j^{1/2}a_\mu$ with $\lambda_j a_\mu a_\nu \leq r_{\mu\nu}$ for all indices $\mu<\nu$ of block $I_j$. Finally, suitable numbers $\vartheta_{ij} \in (0,1)$ are determined. (Obviously, this provides better inequalities than the trivial choice with identical values $a_\mu a_\nu = \min r_{\mu\nu}$.)

Generally, it is difficult to find for a given $R$ with positive correlations an M-matrix $R_0^{-1}$ with $R_0 < R$ in such a way that good positive and numerically available lower bounds are obtained for the difference $P_R(\bigcap_{i=1}^{n} A_i) - P_{R_0}(\bigcap_{i=1}^{n} A_i)$. Theorem 3 and theorem 5 in conjunction with theorem 4 give some examples for this goal.

In particular, let $R_\Theta$ be a correlation matrix as in (5.4) with $p=3$ and $R > R_\Theta$ any correlation matrix with the same three one-factorial diagonal blocks $R_{ii}$ as in $R_\Theta$. If $R_0$ is obtained from $R_\Theta$ by replacing the three num-bers $\vartheta_i = \vartheta_{jk}$ by zeros then $R > R_\Theta > R_0$ and

$$P_{R_0}(\bigcap_{i=1}^{n} A_i) = P(\bigcap_{i\le n_1} A_i)P(\bigcap_{i=n_1+1}^{n_1+n_2} A_i)P(\bigcap_{i>n_1+n_2} A_i).$$

Then we obtain for the “excess” the inequalities

$$P_R(\bigcap_{i=1}^{n} A_i) - P_{R_0}(\bigcap_{i=1}^{n} A_i) > P_{R_\Theta}(\bigcap_{i=1}^{n} A_i) - P_{R_0}(\bigcap_{i=1}^{n} A_i) > 0, \tag{5.14}$$

where the 2nd excess is given by the series (5.10) (or (5.12)) without the first term with $k=0,$ under the additional condition from (5.7).

E.g. let be given an $(9\times 9)$ - correlation matrix. If we get after a transformation with a suitable signature matrix a correlation matrix with only positive correlations and - possibly after a suitable permutation of the indices - a correlation matrix $R=(R_{jk})$ with nine $(3\times 3)$ - blocks $R_{jk}$ and one-factorial diagonal blocks $(R_{jj})$ with factors $a_\mu \in (0,1),$ then we have to find the correlations $\vartheta_{jk} \in (0,1),\ 1\le j<k\le 3,$ which satisfy the conditions

$$\vartheta_{jk} \le \min\left\{ \tfrac{r_{\mu\nu}}{a_\mu a_\nu} \mid \mu\in I_j, \nu\in I_k \right\},\ I_1=\{1,2,3,\},\ I_2=\{4,5,6,\},\ I_3=\{7,8,9,\}, \tag{5.15}$$

and the condition (5.11) to guarantee an M-matrix $R_\Theta^{-1}$. Then the inequality (5.14) holds. It is presumably useful to bring the higher correlations of the given correlation matrix by a suitable permutation into the diagonal blocks as far as possible.

If an $(n\times n)$ - correlation matrix with only positive correlations is written as a $(2\times 2)$ - block matrix $(R_{ik})$ with $(n_i\times n_i)$ - matrices $R_{ii}$ and if $R_\rho$ is obtained from $R$ by replacing all the correlations within $R_{12}$ and $R_{21}$ by their minimum $\rho$ (or by a smaller value $\rho$) then we can obtain further inequalities for the excess in (1.1) if $R_\rho^{-1}$ is an M-matrix, but actual computations require much effort even for moderate dimensions (apart from simulations). If the $R_{ii}^{-1}$ are diagonal dominant M-matrices with only negative off-diagonal elements then we obtain with $R_{12}=\rho\vec{1}_1\vec{1}_2'$ and $q_i := \vec{1}_i' R_{ii}^{-1}\vec{1}_i$ the inverse matrix

$$R_\rho^{-1} = R_{11}^{-1}\oplus R_{22}^{-1} + \frac{1}{1-\rho^2 q_1 q_2}\begin{pmatrix} \rho^2 q_2 R_{11}^{-1}\vec{1}_1\vec{1}_1' R_{11}^{-1} & -\rho R_{11}^{-1}\vec{1}_1\vec{1}_2' R_{22}^{-1} \\ -\rho R_{22}^{-1}\vec{1}_2\vec{1}_1' R_{11}^{-1} & \rho^2 q_1 R_{22}^{-1}\vec{1}_2\vec{1}_2' R_{22}^{-1} \end{pmatrix}, \tag{5.16}$$

which is an M-matrix for sufficiently small values $\rho.$ With the index subsets $\varnothing \ne J_{1\kappa} \subseteq I_1=\{1,\ldots,n_1\},\ \varnothing\ne J_{2\kappa}$ $\subseteq I_2=\{n_1+1,\ldots,n\}$, the submatrices $R_{J_{i\kappa}}$ with the row and column indices from $J_{i\kappa}$ and

$$q_{i\kappa} = \vec{1}_{J_{i\kappa}}' R_{J_{i\kappa}}^{-1}\vec{1}_{J_{i\kappa}} \left|R_{J_{i\kappa}}\right|$$

we find

$$\begin{aligned} &\left|I+R_\rho T\right| \\ &= 1+\sum_\kappa \left|R_{J_{1\kappa}}\right|\left|T_{J_{1\kappa}}\right| + \sum_\lambda \left|R_{J_{2\lambda}}\right|\left|T_{J_{2\lambda}}\right| + \sum_\kappa\sum_\lambda \left|R_{J_{1\kappa}}\right|\left|R_{J_{2\lambda}}\right|\left|T_{J_{1\kappa}}\right|\left|T_{J_{2\lambda}}\right| - \rho^2\sum_\kappa\sum_\lambda q_{1\kappa}q_{2\lambda}\left|T_{J_{1\kappa}}\right|\left|T_{J_{2\lambda}}\right| \\ &= \left|I_1+R_{11}T_1\right|\left|I_2+R_{22}T_2\right|\left(1-\rho^2\frac{\sum_\kappa q_{1\kappa}\left|T_{J_{1\kappa}}\right|}{\left|I_1+R_{11}T_1\right|}\frac{\sum_\lambda q_{2\lambda}\left|T_{J_{2\lambda}}\right|}{\left|I_2+R_{22}T_2\right|}\right). \end{aligned}$$

Then, with the binomial series for the Lt $\left|I+R_\rho T\right|^{-\alpha}\left|T\right|^{-1}$ of the $\Gamma_\alpha(R_\rho)$ - cdf $F(\vec{x};\alpha,R_\rho)$ and the multinomial formula for the powers we obtain by termwise inversion the following series for the $\Gamma_\alpha(R_\rho)$ - cdf:

$$F(\vec{x};\alpha,R_\rho)=\sum_{k=0}^{\infty}\binom{\alpha+k-1}{k}\left[\prod_{i=1}^{2}\left(k!\sum_{(k)}\left(\prod_{\kappa}\frac{q_{i\kappa}^{k_{i\kappa}}}{k_{i\kappa}!}\right)\left(\prod\nolimits_{\mu\in I_i}\partial_\mu^{M_{i\mu}}\right)F(\vec{x}_i;\alpha+k,R_{ii})\right)\right]\rho^{2k}, \qquad (5.17)$$

where $\partial_\mu$ stands for $\frac{\partial}{\partial x_\mu}$ and $M_{i\mu}=\sum_{\mu\in J_{i\kappa}}k_{i\kappa}$. For small values $\rho$ this leads with the indicator functions $e_{J_{i\kappa}}$ of $J_{i\kappa}$ on $I_i$ to the approximation

$$\begin{aligned}F(\vec{x};\alpha,R_\rho)&=F(\vec{x}_1;\alpha,R_{11})F(\vec{x}_2;\alpha,R_{22})+\alpha\left[\prod\nolimits_{i=1}^{2}\left(\sum\nolimits_{\kappa}q_{i\kappa}\left(\prod\nolimits_{\mu\in J_{i\kappa}}\partial_\mu\right)F(\vec{x}_i;\alpha+1,R_{ii})\right)\right]\rho^2\\&+\frac{\alpha(\alpha+1)}{2}\left[\prod\nolimits_{i=1}^{2}\left(\sum\nolimits_{\kappa}\sum\nolimits_{\lambda}q_{i\kappa}q_{i\lambda}\left(\prod\nolimits_{\mu\in I_i}\partial_\mu^{e_{J_{i\kappa}}(\mu)+e_{J_{i\lambda}}(\mu)}\right)F(\vec{x}_i;\alpha+2,R_{ii})\right)\right]\rho^4+O(\rho^6), \qquad (5.18)\end{aligned}$$

where the terms with identical functions $e_{J_{i\kappa}}+e_{J_{i\lambda}}$ on $I_i$ can be collected.

If the M-matrices $R_{ii}^{-1}$ are not diagonal dominant then we can always find some vectors $\vec{a}_i>\vec{0}_i$ with $R_{ii}^{-1}\vec{a}_i>\vec{0}_i$ and $\vec{a}_1\vec{a}_2'\le R_{12}$. Then

$$R_a^{-1}=\begin{pmatrix}R_{11} & \vec{a}_1\vec{a}_2'\\ \vec{a}_2\vec{a}_1' & R_{22}\end{pmatrix}^{-1}$$

is again an M-matrix for sufficiently small values $\|\vec{a}_i\|$, and we obtain an approximation (and a series) for the cdf $F(\vec{x};\alpha,R_a)$ if we replace $\rho q_{i\kappa}$ in (5.18) by $q_{i\kappa}:=\vec{a}_{J_{i\kappa}}'R_{J_{i\kappa}}^{-1}\vec{a}_{J_{i\kappa}}\left|R_{J_{i\kappa}}\right|$. Many 4x4-correlation matrices are 2-factorial. With 2-factorial $R_{ii}$ we can compute the derivatives of $F(\vec{x}_i;\alpha+k,R_{ii})$ by means of formula (2.9) with $m=2$ and $\frac{\partial}{\partial x}g_{\alpha+k}(x;y)=g_{\alpha+k-1}(x;y)-g_{\alpha+k}(x;y)$. Therefore, formula (5.18) is possibly useful for $n\le 8$.

For low dimensions (e.g. $n\le 10$) and small "p-values" $p=P\{\max X_i>x\}=1-P(\bigcap_{i=1}^{n}A_i)$ the following approximation was recommended in [19]:

$$P(\bigcap\nolimits_{i=1}^{n}A_i)\approx P(\bigcap\nolimits_{i\le n_1}A_i)P(\bigcap\nolimits_{i>n_1}A_i)+\sum\nolimits_{k=0}^{\infty}\binom{\alpha+k-1}{k}\left(\frac{r^2}{r_1r_2}\right)^k c_k(x;\alpha,n_1,r_1)c_k(x;\alpha,n_2,r_2) \qquad (5.19)$$

with $c_k(x;\alpha,n_i,r_i)=\binom{\alpha+k-1}{k}^{-1}\int_0^\infty\left(G_\alpha(\frac{x}{1-r_i};\frac{r_iy}{1-r_i})\right)^{n_i}L_k^{(\alpha-1)}(y)g_\alpha(y)dy$ and positive mean correlations $r_i$ from the diagonal blocks $R_{ii}$ on condition that the mean square correlation $r^2$ from the block $R_{12}$ is not larger than $r_1r_2$. For the computation of a $\Gamma_n(\alpha,R)$- cdf with $n\le 4$ and very accurate approximations for $n=5$ see the appendix in [5]. If $n_1>5$ and (or) $n_2>5$ then $P(\bigcap_{i\le n_1}A_i)$, $P(\bigcap_{i>n_1}A_i)$ can be approximated in a similar way under the corresponding assumptions. Numerical examples seem to show that the absolute error of these (frequently conservative) approximations is often smaller than the deviations with refined Bonferroni inequalities of third or even fourth order. The latter ones are inequalities of the form

$$P(\bigcup\nolimits_{i=1}^{n}\overline{A}_i)<P(\overline{A}_1)+P(\overline{A}_2A_1)+P(\overline{A}_3A_2A_1)+P(\overline{A}_4A_3A_2A_1)+\ldots+P(\overline{A}_nA_{n-1}A_{n-2}A_{n-3})$$

with the complemens $\overline{A}_i$ of $A_i$, which can be frequently improved by suitable permutations of the indices.

In theorem 3 in [19] the following only local inequality for $\Gamma_n(\alpha,R)$- cdfs with identical values $x_i=x$ was proved. We consider the set of all $(n\times n)$- correlation matrices $R$ with the same mean correlation $r\ge 0$. Then for any fixed values $x>0,\ n\ge 3,$ and any admissible $\alpha\ge\frac{1}{2}$ the function $\psi:R\to P\{\max X_i\le x;\alpha,R\}$ has a local minimum within this set at the correlation matrix $R_r$ with identical correlations $r$ if

$\lambda = a+(n-4)b-(n-3)c>0,$ where

$$a=\int_0^\infty(\alpha f_1^2-2ryf_1f_2+2r^2y^2f_2^2)F^{n-2}g_\alpha(y)dy,\ b=r\int_0^\infty f_1^2(2ryf_2-f_1)F^{n-3}yg_\alpha(y)dy,$$

$$c=2r^2\int_0^\infty f_1^4F^{n-4}y^2g_\alpha(y)dy,\ \text{and}\ F:=G_\alpha(\tfrac{x}{1-r};\tfrac{ry}{1-r}),\ f_k=\frac{\partial^k}{\partial x^k}G_{\alpha+k}(\tfrac{x}{1-r};\tfrac{ry}{1-r}),\ k=1,2. \qquad (5.20)$$

Frequently, $\lambda>0$ is recognized already by a plot of the integrand of $\lambda$. This result was obtained by an analysis of the Taylor polynomial of 2nd degree for $\psi(R_r+H)$ with small "correlation deviations" $h_{ij}$ in $H$ and $\sum_{i<j}h_{ij}=0.$

It seems that a corresponding global inequality is hardly available, but further investigations are recommended.

## Appendix

### On the convergence of the series representation of the char. function in (5.13)

The char function is

$$|I_n-iR_\Theta T|^{-\alpha}=\left(\prod_{\mu=1}^n z_\mu^\alpha\right)\left(\prod_{i=1}^3\delta_i^{-\alpha}\right)\left(1-(\vartheta_{12}^2\zeta_1\zeta_2+\vartheta_{13}^2\zeta_1\zeta_3+\vartheta_{23}^2\zeta_2\zeta_3-2\vartheta_{12}\vartheta_{13}\vartheta_{23}\zeta_1\zeta_2\zeta_3)\right)^{-\alpha} \qquad (A1)$$

with complex values $\zeta_j=1-\delta_j^{-1}$ and $z_\mu=(1-i(1-a_\mu^2)t_\mu)^{-1}$. The values $z_\mu$ and $u_\mu=1-z_\mu$ lie on the circle defined by $e^{i\gamma}\cos\gamma,\ |\gamma|\le\frac{1}{2}\pi.$ The values $\zeta_j$ lie on the circles given by

$$d_je^{i\gamma}\cos(\gamma),\ d_j\le\left(\sum_{\mu\in I_j}\frac{a_\mu^2}{1-a_\mu^2}\right)\left(1+\sum_{\mu\in I_j}\frac{a_\mu^2}{1-a_\mu^2}\right)^{-1},\ j=1,2,3. \qquad (A2)$$

For the convergence of the binomial series for the last factor in (A1) we need with $\tilde\vartheta_i=\vartheta_{jk}\sqrt{d_jd_k}$ and $c_j=\cos\gamma_j,\ s_j=\sin\gamma_j,\ t_j=\tan\gamma_j$ the condition

$$\rho^2=$$

$$|\vartheta_{12}^2d_1d_2c_1c_2e^{i(\gamma_1+\gamma_2)}+\vartheta_{13}^2d_1d_3c_1c_3e^{i(\gamma_1+\gamma_3)}+\vartheta_{23}^2d_2d_3c_2c_3e^{i(\gamma_2+\gamma_3)}-2\vartheta_{12}\vartheta_{13}\vartheta_{23}d_1d_2d_3c_1c_2c_3e^{i(\gamma_1+\gamma_2+\gamma_3)}|^2=$$

$$|\tilde\vartheta_3^2c_1c_2e^{-i\gamma_3}+\tilde\vartheta_2^2c_1c_3e^{-i\gamma_2}+\tilde\vartheta_1^2c_2c_3e^{-i\gamma_1}-2\tilde\vartheta_1\tilde\vartheta_2\tilde\vartheta_3c_1c_2c_3|^2=$$

$$\frac{\tilde\vartheta^2+(\sum_{j=1}^3\tilde\vartheta_j^2t_j)^2}{\prod_{j=1}^3(1+t_j^2)}<1,\ \text{where}\ \tilde\vartheta=1-|\tilde\Theta|=\sum_{j=1}^3\tilde\vartheta_j^2-2\tilde\vartheta_1\tilde\vartheta_2\tilde\vartheta_3\in(0,1). \qquad (A3)$$

This inequality is not always satisfied, in particular for larger values of the correlations $\tilde\vartheta_j$. The difference between the denominator and the nominator of $\rho^2$ is

$$1-\tilde\theta^2+\sum_{j=1}^3(1-\tilde\vartheta_j^4)t_j^2-2\sum_{i<j}\tilde\vartheta_i^2\tilde\vartheta_j^2t_it_j+\sum_{i<j}t_i^2t_j^2+(t_1t_2t_3)^2.$$

The here appearing quadratic form is positive semi-definite iff

$$\sum_{j=1}^3\tilde\vartheta_j^4\le1. \qquad (A4)$$

Therefore, all $\rho<1$ is guaranteed under this condition with the upper bounds for the $d_j$ from (A2).